# DYNAMICAL MODELS FOR CIRCLE COVERING: BROWNIAN MOTION AND POISSON UPDATING

BY JOHAN JONASSON[1] AND JEFFREY E. STEIF[1,2]

*Chalmers University of Technology and Göteborg University*

We consider two dynamical variants of Dvoretzky's classical problem of random interval coverings of the unit circle, the latter having been completely solved by L. Shepp. In the first model, the centers of the intervals perform independent Brownian motions and in the second model, the positions of the intervals are updated according to independent Poisson processes where an interval of length $\ell$ is updated at rate $\ell^{-\alpha}$ where $\alpha \geq 0$ is a parameter. For the model with Brownian motions, a special case of our results is that if the length of the $n$th interval is $c/n$, then there are times at which a fixed point is not covered if and only if $c < 2$ and there are times at which the circle is not fully covered if and only if $c < 3$. For the Poisson updating model, we obtain analogous results with $c < \alpha$ and $c < \alpha + 1$ instead. We also compute the Hausdorff dimension of the set of exceptional times for some of these questions.

## 1. Introduction.

1.1. *The classical* (*static*) *circle covering model.* Let $C$ denote the circle with circumference 1 and consider a decreasing sequence $\{\ell_n\}_{n \geq 1}$ of positive numbers approaching 0. Let $\{U_n\}_{n \geq 1}$ be a sequence of independent random variables each of which is uniformly distributed on $C$. Let $I_n$ be the open arc of $C$ with center point $U_n$ and length $\ell_n$. Let $E := \limsup_n I_n$ and $F := E^c$. It follows immediately from the Borel–Cantelli lemma that for each $x \in C$, $P(x \in E) = 1$ if and only if $\sum_{n=1}^{\infty} \ell_n = \infty$. Fubini's theorem yields immediately that in this case $F$ has Lebesgue measure 0 a.s. In 1956, Dvoretzky

Received January 2007; revised January 2007.
[1]Supported in part by the Swedish Research Council.
[2]Supported in part by the Göran Gustafsson Foundation for Research in Natural Sciences and Medicine.

(see [4]) raised the question of whether in the $\sum_n \ell_n = \infty$ case it was possible that $F$ was nonempty and gave examples where this occurred. There were a number of various contributions to this question with the final result proved by Shepp (see [16]). Note that Kolmogorov's 0–1 law tells us that $P(F = \varnothing) \in \{0, 1\}$.

THEOREM 1.1 [16].  $P(F = \varnothing) = 1$ *if and only if*

$$\sum_{n=1}^{\infty} \frac{1}{n^2} e^{\ell_1 + \cdots + \ell_n} = \infty.$$

*In particular, if $\ell_n = c/n$ for all $n$, then $P(F = \varnothing) = 1$ if and only if $c \geq 1$.*

The special cases $\ell_n = c/n$ for a constant $c$ were known earlier. The result for $c > 1$ was proved by Kahane (see [10]) and that for $c < 1$ was proved by Billard (see [3]). For the case $c = 1$, Billard also showed that $F$ is at most countable while Mandelbrot (see [13]) and independently Orey (unpublished) then showed that $F$ is a.s. empty in this case. The result that $F$ is at most countable for $c = 1$ also appeared in the first edition of Kahane's book (see [11]) where some of the above results were also presented. The second edition of this book also contains some more history as well as other results such as the Hausdorff dimension of $F$ and a determination of which sets intersect $F$ with positive probability, described in terms of their Hausdorff dimension. We finally mention that in recent years, many refinements of these results have been obtained; see [1, 5, 6]. We finally mention that it is trivial to check that for any sequence $\{\ell_n\}_{n \geq 1}$, $E$ is dense a.s.

1.2. *The dynamical circle covering model.* In this paper, we consider two dynamical variants of the above problem. In the first of these models, each of the centers $U_n$ perform independent Brownian motions on $C$, each with variance 1. In the second model, we associate independent Poisson processes with the different intervals, where the Poisson process associated with the $n$th interval has intensity $\ell_n^{-\alpha}$ for some parameter $\alpha \geq 0$. At the times of the Poisson process associated to the $n$th interval, $I_n$ is given a new center, chosen uniformly on $C$, independent of everything else.

We then ask for each of these two models if there are exceptional times at which we see different "covering behavior" from that which is seen in the earlier static model. We have potentially five (or even more) different types of exceptional times, depending on the $\ell_n$'s and which of the two models we are looking at:

(I) times when a fixed point is not covered even though $\sum_n \ell_n = \infty$,
(II) times when the circle is not fully covered even though $\sum_n e^{\ell_1 + \ell_2 + \cdots + \ell_n} / n^2 = \infty$,

(III) times when a fixed point is covered i.o. even though $\sum_n \ell_n < \infty$,
(IV) times when the circle is fully covered i.o. even though $\sum_n e^{\ell_1+\ell_2+\cdots+\ell_n}/n^2 < \infty$,
(V) times when $E$ is not dense.

To state things more formally, consider the first dynamical model. Here we let, for each $i \geq 1$, $\{U_{i,t}\}_{t\geq 0}$ be an independent standard Brownian motion on $C$ started uniformly. For the second dynamical model, let $\{\{C_{i,j}\}_{i,j\geq 1}, \{Y_{i,j}\}_{i,j\geq 1}\}$ be independent random variables with $C_{i,j}$ being uniformly distributed on $C$ and with $Y_{i,j}$ being exponentially distributed with parameter $\ell_i^{-\alpha}$, $\alpha \geq 0$. For each $i \geq 1$, let $T_{i,0} := 0$ and for $j \geq 1$, let $T_{i,j} := \sum_{k=1}^{j} Y_{i,k}$. In this way, $(T_{i,j})_{j\geq 1}$ is a Poisson process with rate $\ell_i^{-\alpha}$, independent for different $i$. Finally, we let $U_{i,t} = \sum_{j=1}^{\infty} C_{i,j} \mathbf{1}_{[T_{i,j-1}, T_{i,j})}(t)$. Henceforth we refer to the first model as the Brownian model and the second as the Poisson model with parameter $\alpha$. In either case, we let $I_{n,t}$ be the open arc of $C$ with center point $U_{n,t}$ and length $\ell_n$. Let $E_t := \limsup_n I_{n,t}$ and $F_t := E_t^c$.

*Motivation.* Dynamical versions of other probabilistic models have previously been studied. Dynamical percolation was initiated in [9] where the edges in ordinary percolation undergo "Poisson updating." In [2], a dynamical version of the Boolean model in continuum percolation was introduced where the centers of the balls undergo independent Brownian motions. The notion of "exceptional times" appears in many other contexts as well, such as the notion of fast and slow points for Brownian motion.

*Conventions.* Our circle $C$ is $\{(x,y) : x^2 + y^2 = 1/(2\pi)^2\}$. When we subtract two elements in $C$, we mean modular arithmetic so that $(1/(2\pi), 0)$ is the identity. If $x \in C$, by $|x|$ we mean arclength from the identity; in this way $|x| \in [0, 1/2]$ and $|x| = 0$ only for $(1/(2\pi), 0)$. The real line projects onto $C$ via $u \to 1/(2\pi)(\cos(2\pi u), \sin(2\pi u))$. We will assume without loss of generality that $\ell_1 \leq 1/2$. Throughout much of the paper, we will also assume that

$$\ell_n = \Theta(1/n), \tag{1}$$

that is, that there are constants $0 < M_0 \leq M_1 < \infty$ such that for every $n$, $M_0/n \leq \ell_n \leq M_1/n$. Besides $\Theta$ notation, as usual $O(1)$ will denote a quantity which is bounded away from $\infty$. In addition, throughout the paper we also put

$$u_n := \prod_{k=1}^{n}(1-\ell_k)$$

and

$$\beta_0 = \inf\left\{\beta : \sum_{n=1}^{\infty} \frac{e^{\ell_1+\cdots+\ell_n}}{n^{1+\beta}} < \infty\right\}.$$

REMARK ON THE PARAMETRIZATION OF THE POISSON MODEL. One might think that the most natural parameter would be $\alpha = 0$. Interestingly, it turns out that the behavior (and results) for the $\alpha = 2$ case matches very well the behavior for the Brownian model. This is due to the fact that the time it takes a Brownian motion to move a distance $\ell_n$ is of order $\ell_n^2$. We therefore thought it was natural to carry out our analysis for general $\alpha \geq 0$. We only consider $\alpha \geq 0$, as $\alpha < 0$ is easily handled and does not lead to any interesting results.

MEASURABILITY REMARK. Insuring the measurability of the events described below can be handled in the same way as was done in [9] for dynamical percolation. Also, the fact that all the events described below have probability 0 or 1 (once we know that they are measurable) follows immediately from Kolmogorov's 0–1 law.

For a fixed point $x \in C$, in the $\ell_n = \Theta(1/n)$ case, it follows immediately from the Borel–Cantelli lemma that $P(x \in F) = 0$ and hence for any of the dynamical models, by Fubini's theorem, $\{t : x \in F_t\}$ has Lebesgue measure 0 a.s. The question we address in the first two theorems is when there are exceptional times $t$ at which $x$ is covered by only finitely many of the $I_{n,t}$'s; that is, $x \in F_t$. See [11] for the definition of Hausdorff dimension which we denote here by HD.

THEOREM 1.2. *Assume that* (1) *holds. Consider the Brownian model and fix $x \in C$.*

(i) *If $\liminf_n n^2 u_n < \infty$, then $P(\exists t \in [0,1] : x \in F_t) = 0$. In particular, if $\ell_n = c/n$ for all $n$, then this holds if $c \geq 2$.*

(ii) *If $\sum_{n=1}^{\infty} e^{\ell_1 + \ell_2 + \cdots + \ell_n}/n^3 < \infty$, then $P(\exists t \in [0,1] : x \in F_t) = 1$. In particular, if $\ell_n = c/n$ for all $n$, then this holds if $c < 2$.*

(iii) *We have that*

$$\mathrm{HD}(\{t \in [0,1] : x \in F_t\}) = \left(1 - \frac{\beta_0}{2}\right) \wedge 0 \qquad a.s.$$

*In particular, in the case $\ell_n = c/n$ for all $n$ with $c \leq 2$, we have*

$$\mathrm{HD}(\{t \in [0,1] : x \in F_t\}) = 1 - \frac{c}{2} \qquad a.s.$$

REMARK. Unfortunately we have not been able to determine the behavior of the Brownian model for the "intermediate" cases when the conditions in (i) and (ii) both fail. An example of such a sequence would be $\ell_n = 2/n - 1/(n \log n)$. On the other hand, an example of a sequence which leads to exceptional times but where the HD of these exceptional times is 0 is given by $\ell_n = 2/n - 1/(n\sqrt{\log n})$.

The Poisson model, however, turns out to be more amenable to our analysis and we obtain an exact condition for having exceptional times of type (I).

THEOREM 1.3. *Assume that* (1) *holds. Consider the Poisson model with parameter $\alpha > 0$. Fix $x \in C$.*

(i) *Then $P(\exists t \in [0,1] : x \in F_t) = 1$ if and only if*
$$\sum_{n=1}^{\infty} \frac{e^{\ell_1 + \ell_2 + \cdots + \ell_n}}{n^{1+\alpha}} < \infty.$$

*In particular, if $\ell_n = c/n$ for all $n$, then this holds if and only if $c < \alpha$.*

(ii) *We have that*
$$\mathrm{HD}(\{t \in [0,1] : x \in F_t\}) = \left(1 - \frac{\beta_0}{\alpha}\right) \wedge 0 \qquad a.s.$$

*In particular, in the case $\ell_n = c/n$ for all $n$ with $c \leq \alpha$, we have*
$$\mathrm{HD}(\{t \in [0,1] : x \in F_t\}) = 1 - \frac{c}{\alpha} \qquad a.s.$$

REMARK. The case $\alpha = 0$ is almost trivially covered by the Borel–Cantelli lemma by noting that the probability that the $n$th interval covers $x$ for the whole time span $[0,1]$ is then at least $e^{-1} \ell_n$. Hence there are no exceptional times of type (I) for $\alpha = 0$.

Our next two results deal with the question of exceptional times of type (II).

THEOREM 1.4. *Assume that* (1) *holds and consider the Brownian model.*

(i) *If $\liminf_n n^3 u_n < \infty$, then*
$$P(\exists t \in [0,1] : F_t \neq \varnothing) = 0.$$

*In particular, if $\ell_n = c/n$ for all $n$, then this holds if $c \geq 3$.*

(ii) *If $\sum_{n=1}^{\infty} e^{\ell_1 + \ell_2 + \cdots + \ell_n} / n^4 < \infty$, then*
$$P(\exists t \in [0,1] : F_t \neq \varnothing) = 1.$$

*In particular, if $\ell_n = c/n$ for all $n$, then this holds if $c < 3$.*

(iii) *We have a.s.:*

(a)
$$\mathrm{HD}(\{(t,x) : x \in F_t\}) = \begin{cases} 2 - \dfrac{\beta_0}{2}, & \text{if } 0 \leq \beta_0 \leq 2, \\ 3 - \beta_0, & \text{if } 2 \leq \beta_0 \leq 3, \\ 0, & \text{if } \beta_0 \geq 3, \end{cases}$$

(b)
$$\mathrm{HD}(\{x : \exists t : x \in F_t\}) \begin{cases} = 1, & \text{if } 0 \leq \beta_0 < 2, \\ \leq 3 - \beta_0, & \text{if } 2 \leq \beta_0 \leq 3, \\ = 0, & \text{if } \beta_0 \geq 3, \end{cases}$$

(c)
$$\mathrm{HD}(\{t : F_t \neq \varnothing\}) \begin{cases} = 1, & \text{if } 0 \leq \beta_0 < 1, \\ \leq \dfrac{3 - \beta_0}{2}, & \text{if } 1 \leq \beta_0 \leq 3, \\ = 0, & \text{if } \beta_0 \geq 3. \end{cases}$$

*In particular, in the case $\ell_n = c/n$ for all $n$ and $c < 3$, then the dimension bounds are simply obtained by plugging in $c$ for $\beta_0$.*

REMARK. The first equalities in (b) and (c) hold since the event in question then occurs at a fixed time. Note the lack of smoothness in (a) at $\beta_0 = 2$ which is of course due to the fact that 2 is the critical value arising in Theorem 1.2. As for the type (I) case, there are intermediate cases such as $\ell_n = 3/n - 1/(n \log n)$ where both (i) and (ii) fail and so we cannot determine if there are exceptional times. This will also occur in the Poisson case.

THEOREM 1.5. *Assume that* (1) *holds. Consider the Poisson model with parameter $\alpha > 0$.*

(i) *If* $\liminf_n n^{1+\alpha} u_n < \infty$, *then*
$$P(\exists t \in [0, 1] : F_t \neq \varnothing) = 0.$$

*In particular, if $\ell_n = c/n$ for all $n$, then this holds if $c \geq 1 + \alpha$.*

(ii) *If*
$$\sum_{n=1}^{\infty} \frac{e^{\ell_1 + \ell_2 + \cdots + \ell_n}}{n^{2+\alpha}} < \infty,$$

*then $P(\exists t \in [0, 1] : F_t \neq \varnothing) = 1$. In particular, when $\ell_n = c/n$ for all $n$, then this holds if $c < 1 + \alpha$.*

(iii) *We have a.s.:*

(a) *(for $\alpha \geq 1$)*
$$\mathrm{HD}(\{(t, x) : x \in F_t\}) = \begin{cases} 2 - \dfrac{\beta_0}{\alpha}, & \text{if } 0 \leq \beta_0 \leq \alpha, \\ 1 + \alpha - \beta_0, & \text{if } \alpha \leq \beta_0 \leq 1 + \alpha, \\ 0, & \text{if } \beta_0 \geq 1 + \alpha, \end{cases}$$

(a′) [*for* $\alpha \in (0,1)$]

$$\mathrm{HD}(\{(t,x):x\in F_t\}) = \begin{cases} 2-\beta_0, & \text{if } 0\leq \beta_0 \leq 1, \\ \dfrac{1+\alpha-\beta_0}{\alpha}, & \text{if } 1\leq \beta_0 \leq 1+\alpha, \\ 0, & \text{if } \beta_0 \geq 1+\alpha, \end{cases}$$

(b)

$$\mathrm{HD}(\{x:\exists t:x\in F_t\}) \begin{cases} =1, & \text{if } 0\leq \beta_0 < \alpha, \\ \leq 1+\alpha-\beta_0, & \text{if } \alpha \leq \beta_0 \leq 1+\alpha, \\ =0, & \text{if } \beta_0 \geq 1+\alpha, \end{cases}$$

(c)

$$\mathrm{HD}(\{t:F_t\neq \varnothing\}) \begin{cases} =1, & \text{if } 0\leq \beta_0 < 1, \\ \leq \dfrac{1+\alpha-\beta_0}{\alpha}, & \text{if } 1\leq \beta_0 \leq 1+\alpha, \\ =0, & \text{if } \beta_0 \geq 1+\alpha. \end{cases}$$

*In particular, in the case $\ell_n = c/n$ for all $n$ and $c < 1+\alpha$, then the dimension bounds are simply obtained by plugging in $c$ for $\beta_0$.*

REMARK. The difference in the form of the Hausdorff dimension in (a) and (a′) is due to the fact that as $\beta_0$ decreases starting from $\infty$, when $\alpha > 1$, we encounter exceptional points on the circle in the sense of Theorem 1.2 before we encounter exceptional times in $[0,1]$ in the sense of Theorem 1.1, while when $\alpha < 1$, we encounter these objects in the opposite order. As in Theorem 1.4, there are intermediate cases such as $\ell_n = (\alpha+1)/n - 1/(n\log n)$ where both (i) and (ii) fail and so we cannot determine if there are exceptional times.

As for type (I) exceptional times, the case $\alpha = 0$ requires special treatment, but unlike the type (I) situation, it is not trivial. Indeed there are situations with $\alpha = 0$ where the circle is fully covered i.o. in the static model, but where there are exceptional times at which some point on the circle fails to be covered infinitely often; the sufficient condition differs from Shepp's condition for the static case by a logarithmic factor.

THEOREM 1.6. *Assume that* (1) *holds and consider the Poisson model with $\alpha = 0$.*

(i) *If $\liminf_n n(\log n) u_n < \infty$, then $P(\exists t \in [0,1]: F_t \neq \varnothing) = 0$.*
(ii) *If*

$$\sum_{n=1}^{\infty} \frac{e^{\ell_1+\ell_2+\cdots+\ell_n}}{n^2 \log n} < \infty,$$

*then $P(\exists t \in [0,1]: F_t \neq \varnothing) = 1$.*

REMARK. An example of a sequence where $F = \varnothing$ in the static model but for which there are exceptional times from this is given by $\ell_n = 1/n - 1/(n \log n)$. Note that the case $\ell_n = 1/n$ is not covered by parts (i) or (ii) and so we cannot determine if there are exceptional times in this case. We mention that one can also prove, along the same lines as the other HD results, that $\mathrm{HD}(\{(t,x) : x \in F_t\}) \leq 1$; this bound is also strongly suggested by Theorem 1.5(iii)(a').

We now move to results concerning the $\sum_{n=1}^{\infty} \ell_n < \infty$ case, which we feel are less central than the $\sum_{n=1}^{\infty} \ell_n = \infty$ results. We start with type (III) exceptional times.

THEOREM 1.7. *Assume that $\sum_n \ell_n < \infty$, fix $x \in C$ and let $T := \{t \in [0,1] : x \in E_t\}$.*

(a) *In the Brownian model, $P(T \neq \varnothing) = 1$.*
(b) *In the Poisson model, if $\sum_n \ell_n^{1-\alpha} < \infty$, then $P(T = \varnothing) = 1$ while if $\sum_n \ell_n^{1-\alpha} = \infty$, then $P(T \neq \varnothing) = 1$.*

For type (IV) exceptional times, we have no results but finally for type (V), we have the following.

THEOREM 1.8. *Let $T := \{t \in [0,1] : E_t$ is not dense$\}$.*

(a) *For the Brownian model and any sequence $\{\ell_n\}_{n \geq 1}$, we have that $P(T \neq \varnothing) = 0$.*
(b) *In the Poisson model with $\alpha > 0$, if $\ell_n \geq 1/n^c$ for all $n$ and some $c$, then $P(T \neq \varnothing) = 0$. For $\alpha = 0$, we have that $P(T \neq \varnothing) = 0$ for all $\{\ell_n\}_{n \geq 1}$.*
(c) *In the Poisson model with $\alpha > 0$, there exists a sequence $\{\ell_n\}_{n \geq 1}$ so that $P(T \neq \varnothing) = 1$.*

REMARK. For our results, in obtaining both existence of exceptional times and *lower* bounds on Hausdorff dimension, the key step is to obtain bounds on various correlations. For these parts of Theorems 1.2 and 1.3, once one has good bounds on correlations, one could place things in the context of Proposition A.16 in [12] (or as is done in [15]). Similarly, Proposition A.13 in [12], on the other hand, can be used, once certain bounds are obtained, to give upper bounds on Hausdorff dimension (or nonexistence of certain exceptional times). Most of our arguments, however, will not explicitly put things in that context, although we will skip certain Hausdorff dimension arguments and simply state that they can be put into this context. For Theorems 1.4, 1.5 and 1.6, we have, however, space–time results which do not fit as well into the context of Proposition A.16 since the latter is stated for "one-dimensional" situations. In our case, we also need to treat the time and space coordinates differently. Nonetheless, one should think that the ideas of these Propositions A.13 and A.16 are always lurking in the background.

The rest of the paper is organized as follows. In Section 2, we prove Theorems 1.2 and 1.3 and in Section 3, we prove Theorems 1.4, 1.5 and 1.6. These correspond respectively to type (I) and type (II) results. The proofs of type (III) and type (V) results are fairly easy and we therefore omit them. For example, for Theorem 1.8 it suffices to observe that at any fixed time and for any of the arcs, there will a.s. be an arc contained in it, which in turn contains another arc, and so on. For more details on this and on some of the other proofs, we refer to an extended version of this article that can be found on either of the authors' homepages.

**2. Proofs of type (I) results.** Recall our standing assumption (1). We begin with three technical lemmas that will prove useful on several occasions.

LEMMA 2.1. *Assume that* (1) *holds. Let* $\beta > 0$. *Then for every* $b > 0$

$$(2) \qquad \int_0^b e^{\sum_{n:\ell_n^\beta \geq t} \ell_n} \, dt < \infty$$

*if and only if*

$$\sum_{n=1}^\infty \frac{e^{\ell_1 + \ell_2 + \cdots + \ell_n}}{n^{1+\beta}} < \infty.$$

PROOF. We will use Lemma 11.4.1 of [11] which states that for a convex decreasing function $f(t)$ on $(0, b)$, $\int_0^b e^{f(t)} \, dt < \infty$ if and only if $\int_0^b e^{f(t)} \frac{df'(t)}{f'(t)^2} < \infty$.

To apply this result, put $f(t) = \sum_{n=1}^\infty \ell_n^{1-\beta}(\ell_n^\beta - t)^+$. Then $f(t)$ is decreasing and convex and

$$\left| f(t) - \sum_{n : \ell_n^\beta \geq t} \ell_n \right| = t \sum_{n : \ell_n^\beta \geq t} \ell_n^{1-\beta} = O(1) t \sum_{k=1}^{(M_1/t^{1/\beta})} k^{\beta - 1} = O(1).$$

Hence (2) is equivalent to $\int_0^b e^{f(t)} \, dt < \infty$. We now use the above result. We have

$$f'(t) = -\sum_{k=1}^n \ell_k^{1-\beta}, \qquad \ell_{n+1}^\beta < t < \ell_n^\beta,$$

and in particular that $f'(t) = -\Theta(n^\beta)$, $\ell_{n+1}^\beta < t < \ell_n^\beta$. Since $f'(t)$ jumps when $t = \ell_n^\beta$ and the size of the corresponding jump is $\ell_n^{1-\beta} = \Theta(n^{\beta-1})$, we get that

$$\int_0^b e^{\sum_{n : \ell_n^\beta \geq t} \ell_n} \, dt < \infty$$

if and only if
$$\sum_{n=1}^{\infty} \frac{e^{f(\ell_n^{\beta})}}{n^{1+\beta}} < \infty.$$

Finally
$$f(\ell_n^{\beta}) = \ell_1 + \ell_2 + \cdots + \ell_n - \ell_n^{\beta} \sum_{k=1}^{n} \ell_k^{1-\beta} = \ell_1 + \ell_2 + \cdots + \ell_n + O(1)$$

and the lemma follows. □

LEMMA 2.2.
$$\sup\left\{\beta : \liminf_n n^{\beta} u_n < \infty\right\} \geq \inf\left\{\beta : \sum_{n=1}^{\infty} \frac{e^{\ell_1+\ell_2+\cdots+\ell_n}}{n^{1+\beta}} < \infty\right\}$$
$$= \limsup_n \frac{\ell_1 + \ell_2 + \cdots + \ell_n}{\log n}.$$

REMARK. The last expression comes up in Section 11.8 in [11] where the HD of the set $F$ is studied. If $\ell_n = \Theta(1/n)$, then Steps 1 and 2 in the proof of Theorem 1.2 tell us that the first inequality is an equality.

PROOF OF LEMMA 2.2. For the first inequality, it suffices to show that for any $\varepsilon > 0$,
$$\liminf_n n^{\beta} u_n = \infty \quad \text{implies that} \quad \sum_{n=1}^{\infty} \frac{e^{\ell_1+\ell_2+\cdots+\ell_n}}{n^{1+\beta+\varepsilon}} < \infty.$$

However, this follows from noting that $e^{-x} \geq 1-x$ implies that $e^{\ell_1+\ell_2+\cdots+\ell_n} \leq 1/u_n$.

Let $L$ denote the third expression. Fix $\varepsilon > 0$. Then $\ell_1 + \ell_2 + \cdots + \ell_n \leq (L+\varepsilon/2)\log n$ and hence $e^{\ell_1+\ell_2+\cdots+\ell_n} \leq n^{L+\varepsilon/2}$ for all $n$ large. It follows that
$$\sum_{n=1}^{\infty} \frac{e^{\ell_1+\ell_2+\cdots+\ell_n}}{n^{1+L+\varepsilon}} < \infty$$

and hence the third term is at least the second term. For the other direction, we may assume that $L > 0$ and we need to show that for all $\varepsilon \in (0, L)$,
$$\sum_{n=1}^{\infty} \frac{e^{\ell_1+\ell_2+\cdots+\ell_n}}{n^{1+L-\varepsilon}} = \infty.$$

We have $\ell_1 + \ell_2 + \cdots + \ell_n \geq (L-\varepsilon)\log n$ and hence $e^{\ell_1+\ell_2+\cdots+\ell_n} \geq n^{L-\varepsilon}$ i.o. It is not hard to show that if $e^{\ell_1+\ell_2+\cdots+\ell_{n_0}} \geq n_0^{L-\varepsilon}$, then
$$\sum_{k=n_0}^{2n_0} \frac{e^{\ell_1+\ell_2+\cdots+\ell_k}}{k^{1+L-\varepsilon}} = \frac{1}{O(1)},$$

the $O(1)$ term being independent of $n_0$. To do this, one simply bounds $e^{\ell_1+\ell_2+\cdots+\ell_k}$ from below by $n_0^{L-\varepsilon}$ for each $k$ and computes. This clearly implies divergence of the series. □

The last lemma is elementary and the proof is left to the reader.

LEMMA 2.3. *There exists a constant $C$ so that for all $t, a \in [0,1]$ and $b \in (-1/2, 1/2)$, if $Z^*$ is a normal random variable with mean $0$ and variance $t$, then*
$$\sum_{k \in \mathbb{Z}} P(Z^* \in (k+b-a, k+b+a)) \leq C P(Z^* \in (b-a, b+a)).$$

PROOF OF THEOREM 1.2. We begin with part (i). Let $I_t = \bigcup_{n=1}^\infty I_{n,t}$ and $J_t = I_t^c$ and note that it is elementary that $P(\exists t \in [0,1] : x \in F_t) = 1$ if and only if $P(\exists t \in [0,1] : x \in J_t) > 0$. Let $T = \{t \in [0,1] : x \in J_t\}$. Now put $J_{n,t} = (\bigcup_{k=1}^n I_{k,t})^c$ and note that $\bigcap_n J_{n,t} = J_t$. We shall first show that if $\liminf_n n^2 u_n = 0$, then

(3) $$\lim_{n \to \infty} P(\exists t : x \in J_{n,t}) = 0.$$

Fix $n \geq 1$. For $i = 1, 2, \ldots, n^2$, let $A_i$ be the event that there exists a $t \in [(i-1)/n^2, i/n^2]$ for which $x \in J_{n,t}$. (We suppress the dependence on $n$ in the notation.) Then
$$P(\exists t : x \in J_{n,t}) = P\left(\bigcup_{i=1}^{n^2} A_i\right) \leq \sum_{i=1}^{n^2} P(A_i) = n^2 P(A_1).$$

For $k \leq n$, let $B_k$ be the event that the $k$th interval covers $x$ for the whole time interval $[0, 1/n^2]$. This event contains the event that $[x-M, x+M] \subseteq I_{k,0}$ where $M = \max_{t \in [0, 1/n^2]} |U_{k,t} - U_{k,0}|$. Letting $B_k'$ denote the latter event, we have
$$P(B_k) \geq P(B_k') = \mathbb{E}[P(B_k'|M)] \geq \mathbb{E}[\ell_k - 2M] = \ell_k - 2\mathbb{E}M \geq \ell_k - \frac{4}{n}$$
where we have used the usual scaling property of Brownian motion. Since $A_1 \subseteq \bigcap_{k=1}^n B_k^c$, the $B_k$'s are independent and $\ell_1 \leq 1/2$, we get
$$P(A_1) \leq \prod_{k=1}^n \left(1 - \ell_k + \frac{4}{n}\right)$$
$$\leq \prod_{k=1}^n (1 - \ell_k)\left(1 + \frac{8}{n}\right) \leq O(u_n).$$

Hence
$$P(\exists t : x \in J_{n,t}) \leq O(n^2 u_n).$$

Since the left-hand side is decreasing in $n$, (3) is established. The case $\liminf_n n^2 u_n \in (0, \infty)$ requires one extra step. Let $N_n$ be the number of $i \in \{1, 2, \ldots, n^2\}$ such that $A_i$ occurs. Then the above arguments show that $\liminf_n \mathbb{E} N_n < \infty$. It is easy to see that $T \leq \liminf_n N_n$. Hence by Fatou's lemma, $\mathbb{E}|T| < \infty$ and so $|T| < \infty$ a.s. Since our process is a reversible stationary Markov process, we finally conclude that $P(T \neq \varnothing) = 0$ by combining Theorem 6.7 in [8] and (2.9) in [7]. This finishes the proof of part (i).

For part (ii), letting $T_n = \{t \in [0,1] : x \notin \bigcup_{k=1}^n I_{n,t}\}$, we have that $T = \bigcap_{n=1}^\infty T_n$. Next, since the intervals are taken to be open and Brownian motion has continuous paths, it follows that the sets $T_n$ are closed and therefore by compactness, $T$ is nonempty if and only if all the $T_n$'s are. Thus if it can be shown that $P(T_n \neq \varnothing)$ is bounded away from 0, then it follows that $P(T \neq \varnothing) > 0$ and (ii) then follows.

Let
$$X_n := \int_0^1 I_{\{t \in T_n\}} \, dt,$$

which is the Lebesgue amount of time that $x$ is not covered by the first $n$ intervals. Since the probability that $x \notin \bigcup_{k=1}^n I_{k,t}$ at a fixed time $t$ is $\prod_{k=1}^n (1 - \ell_k) = u_n$, it follows from Fubini's theorem that $\mathbb{E} X_n = u_n$. We will now establish that $\mathbb{E}[X_n^2] \leq O(u_n^2)$ if (and in fact only if) $\sum_n e^{\ell_1 + \cdots + \ell_n}/n^3 < \infty$. Once this has been done it then follows, under this condition, from the second moment method that
$$P(T_n \neq \varnothing) \geq P(X_n > 0) \geq \frac{(\mathbb{E} X_n)^2}{\mathbb{E}[X_n^2]}$$

is bounded away from 0, as desired. Now, by Fubini's theorem

(4) $$\mathbb{E}[X_n^2] = \int_0^1 \int_0^1 P(\{s \in T_n\} \cap \{t \in T_n\}) \, ds \, dt.$$

By stationarity, it is easy to see that (4) is at most
$$2 \int_0^1 P(\{t \in T_n\} \cap \{0 \in T_n\}) \, dt$$

and at least
$$1/2 \int_0^{1/2} P(\{t \in T_n\} \cap \{0 \in T_n\}) \, dt.$$

Fix $n$. Put $A_t = \{t \in T_n\}$ and $A_{k,t}$ for the event that $x$ is not covered by $I_k$ at time $t$ and note that $A_t = \bigcap_{k=1}^n A_{k,t}$. Clearly
$$P(A_t \cap A_0) = \prod_{k=1}^n P(A_{k,t} \cap A_{k,0}).$$

The probability $P(A_{k,t} \cap A_{k,0})$ is the probability that $I_{k,0} \cap \{x, x - Z_t\} = \varnothing$, where $Z_t := U_{k,t} - U_{k,0}$, the increment of the $k$th interval during the time $[0,t]$. Note that $Z_t$ is a normal random variable with zero mean and variance $t$ projected onto $C$ as described in the Introduction. We have

$$(5) \quad P(A_{k,t} \cap A_{k,0}) = \mathbb{E}[P(A_{k,t} \cap A_{k,0}|Z_t)] = 1 - 2\ell_k + \mathbb{E}[(\ell_k - |Z_t|)^+].$$

Some elementary considerations (using again that $\ell_1 \leq 1/2$) allow us to write (5) as

$$(1-l_k)^2 e^{\mathbb{E}[(\ell_k - |Z_t|)^+]}(1+r_k),$$

where $|r_k| \leq 5\ell_k^2$. We then have that $\prod_{n=1}^{\infty}(1+r_n) < \infty$ and so it follows that $\mathbb{E}[X_n^2] \leq O(u_n^2)$ if

$$(6) \quad \int_0^1 e^{\sum_{n=1}^{\infty} \mathbb{E}[(\ell_n - |Z_t|)^+]} dt < \infty$$

and only if

$$(7) \quad \int_0^{1/2} e^{\sum_{n=1}^{\infty} \mathbb{E}[(\ell_n - |Z_t|)^+]} dt < \infty.$$

Note that trivially

$$\mathbb{E}[(\ell_n - |Z_t|)^+] \leq P(|Z_t| \leq \ell_n)\ell_n$$

and a standard bound on the normal distribution together with Lemma 2.3 gives

$$(8) \quad P(|Z_t| \leq \ell_n) = O(1)\frac{\ell_n}{\sqrt{t}}.$$

These easily yield that (6) holds if and only if (7) does and so we concentrate only on (6).

Next, we have

$$\mathbb{E}[(\ell_n - |Z_t|)^+] = P(|Z_t| \leq \ell_n)(\ell_n - \mathbb{E}[|Z_t| \, | \, 0 \leq |Z_t| \leq \ell_n]).$$

Since a nonnegative random variable conditioned on being smaller than some value is stochastically dominated by the original random variable, we have that the expectation in the right-hand side is bounded above by $\ell_n \wedge \sqrt{t}$. Hence

$$\mathbb{E}[(\ell_n - |Z_t|)^+] \geq P(|Z_t| \leq \ell_n)(\ell_n - \sqrt{t})^+.$$

Using that $\ell_n = \Theta(1/n)$ we get

$$\sum_{n:\ell_n^2 \geq t} \sqrt{t} = O(1).$$

From (8), we get

$$\sum_{n:\ell_n^2<t}\ell_n P(|Z_t|\leq\ell_n)<\frac{1}{\sqrt{t}}\sum_{n:\ell_n^2<t}\ell_n^2=O(1)$$

[where we again used that $\ell_n=\Theta(1/n)$]. Putting this together we have that $\mathbb{E}[X_n^2]\leq O(u_n^2)$ if and only if

$$\int_0^1 e^{\sum_{n:\ell_n^2\geq t}\ell_n P(|Z_t|\leq\ell_n)}\,dt<\infty.$$

Since the probability that a standard normal random variable exceeds a value $y>0$ is bounded above by $O(1)/y$, we get that $\sum_{n:\ell_n^2\geq t}\ell_n P(|Z_t|>\ell_n)\leq\sum_{n:\ell_n^2\geq t}\sqrt{t}=O(1)$ and so the above integral is finite if and only if

$$\int_0^1 e^{\sum_{n:\ell_n^2\geq t}\ell_n}\,dt<\infty.$$

Plugging this into Lemma 2.1 yields $\mathbb{E}[X_n^2]\leq O(u_n^2)$ if and only if

$$\sum_{n=1}^\infty \frac{e^{\ell_1+\ell_2+\cdots+\ell_n}}{n^3}<\infty.$$

This finishes the proof of (ii).

For part (iii), the two key steps are:

*Step* 1. $\liminf_n n^\beta u_n<\infty$ implies that $\mathrm{HD}(T)\leq 1-\beta/2$ a.s.

*Step* 2. $\sum_{n=1}^\infty e^{\ell_1+\ell_2+\cdots+\ell_n}/n^{1+\beta}<\infty$ implies that $P(\mathrm{HD}(T)\geq 1-\beta/2)>0$.

It is clear from Kolmogorov's 0–1 law, the fact that a countable union of sets each of which has HD at most $d$ also has HD at most $d$ and countable additivity that the HD results follow from these two steps and Lemma 2.2.

These steps are really refinements of the arguments already presented. Step 1 follows, after a little work, from what has already been done together with Proposition A.13 in [12] while step 2 follows, after a little work, from what has already been done together with the proof of Proposition A.16 in [12]. We skip the details other than mentioning that one thing which is needed is that Lemma 2.1 tells us that $\sum_{n=1}^\infty e^{\ell_1+\ell_2+\cdots+\ell_n}/n^{1+\beta}<\infty$ implies that

$$\int_0^1 e^{\sum_{n:\ell_n^\beta\geq t}\ell_n}\,dt<\infty$$

which, by an easy change of variables, is equivalent to

(9) $$\int_0^1 e^{\sum_{n:\ell_n^2\geq t}\ell_n}\left(\frac{1}{t}\right)^{1-\beta/2}dt<\infty. \qquad\square$$

In the above proof it was shown that the "Shepp-like" condition $\sum_n e^{\ell_1+\cdots+\ell_n}/n^3 < \infty$ is necessary and sufficient for the second moment argument to work. What we have not been able to determine is if failure of the second moment argument necessarily implies that there are no exceptional times. The reason is that it is difficult to control the conditional distribution of the positions of the first $n$ arcs at the first time when $x$ is not covered by any of them. For the Poisson model this problem vanishes.

PROOF OF THEOREM 1.3. We use exactly the same notation as in the proof of Theorem 1.2. As in that proof, we have $\mathbb{E}X_n = u_n$ and we will show that

$$(10) \qquad E[X_n^2] \leq O(u_n^2) \quad \text{if and only if} \quad \sum_{n=1}^{\infty} e^{\ell_1+\cdots+\ell_n}/n^{1+\alpha} < \infty.$$

However, we first show $E[X_n^2] \leq O(u_n^2)$ is necessary and sufficient for $T$ to be nonempty with positive probability. Note that, using the fact that (14) below is decreasing in $t$, we have

$$(11) \qquad \mathbb{E}[X_n | 0 \in T_n] u_n = \Theta(1) \mathbb{E}[X_n^2].$$

The sufficiency argument is identical to the proof of Theorem 1.2 except for the small irritation that $T_n$ is not a closed set. So $\inf_n P(T_n \neq \varnothing) > 0$ does not immediately allow us to conclude that $P(T \neq \varnothing) > 0$. This very minor issue arose in [9] (as well as in [15]) and was taken care of there by Lemma 3.2. The key observation, left to the reader, is that

$$(12) \qquad \bigcap_{n>0} \overline{T_n} = T \qquad \text{a.s.}$$

This claim takes care of the above problem allowing us to conclude that $P(T \neq \varnothing) > 0$ and will be needed for part (ii) as well.

For the necessity, observe that

$$(13) \qquad P(X_n > 0) = \frac{\mathbb{E}X_n}{\mathbb{E}[X_n | X_n > 0]}.$$

Put $S_n := \min\{t : t \in T_n\}$ (the minimum exists since the Poisson processes are right continuous). Now the crucial observation to make is that at the time $S_n$, the positions of the first $n$ arcs are independent and uniform given that none of them contains $x$. Hence, by conditioning on $S_n$, using translation invariance and the strong Markov property,

$$\mathbb{E}[X_n | X_n > 0] = \Theta(1) \mathbb{E}[X_n | 0 \in T_n].$$

Using this, together with (11) and (13), necessity now follows.

We now show (10). The rest of the proof is very similar to the proof of Theorem 1.2(ii). With $n$ fixed, by conditioning on whether or not arc $I_k$ has been updated before time $t$ we get

(14) $$P(A_{k,t} \cap A_{k,0}) = (1 - \ell_k)^2(1 - e^{-t/\ell_k^\alpha}) + (1 - \ell_k)e^{-t/\ell_k^\alpha}.$$

Hence

$$P(A_t \cap A_0) = u_n^2 \prod_{k=1}^n \left(1 - e^{-t/\ell_k^\alpha} + \frac{e^{-t/\ell_k^\alpha}}{1 - \ell_k}\right)$$

$$= u_n^2 \prod_{k=1}^n \left(1 + \frac{\ell_k e^{-t/\ell_k^\alpha}}{1 - \ell_k}\right)$$

$$= \Theta(1) u_n^2 \prod_{k=1}^n (1 + \ell_k e^{-t/\ell_k^\alpha}),$$

where the $\Theta(1)$ term is bounded between 1 and $\prod_{k=1}^\infty (1 + 2\ell_k^2) < \infty$. Therefore, using (11), we have that $E[X_n^2] \leq O(u_n^2)$ if and only if

$$\int_0^1 \prod_{n=1}^\infty (1 + \ell_n e^{-t/\ell_n^\alpha}) \, dt < \infty.$$

Since $x - x^2 \leq \log(1 + x) \leq x$ on $[0, \infty)$, this is equivalent to

(15) $$\int_0^1 e^{\sum_{n=1}^\infty \ell_n e^{-t/\ell_n^\alpha}} \, dt < \infty.$$

Using $\ell_n \leq M_1/n$ it follows that

$$\sum_{n:\,\ell_n^\alpha \leq t} \ell_n e^{-t/\ell_n^\alpha} \leq t^{1/\alpha} \sum_{n:\,\ell_n^\alpha \leq t} e^{-t/\ell_n^\alpha} \leq t^{1/\alpha} \sum_{n=1}^\infty (e^{-M_1^{-\alpha} t})^{n^\alpha}$$

$$= O(1) t^{1/\alpha} \frac{1}{(1 - e^{-M_1^{-\alpha} t})^{1/\alpha}} = O(1) t^{1/\alpha} \frac{1}{t^{1/\alpha}} = O(1).$$

We also have that

$$\sum_{n:\,\ell_n^\alpha \geq t} \ell_n(1 - e^{-t/\ell_n^\alpha}) \leq t \sum_{n:\,\ell_n^\alpha \geq t} \ell_n^{1-\alpha} = O(1) t \sum_{n=1}^{M_1/t^{1/\alpha}} n^{\alpha - 1} = O(1).$$

Putting this together, (15) is equivalent to

$$\int_0^1 e^{\sum_{n:\,\ell_n^\alpha \geq t} \ell_n} \, dt < \infty.$$

Now Lemma 2.1 finishes the proof of (10).

We skip the proof of part (ii) and simply state that one follows the proof for the Brownian model and that one replaces the $\beta/2$ term by $\beta/\alpha$. $\square$

**3. Proofs of type (II) results.** Recall our standing assumption (1). In this section, the general approach will be to try to analyze the "space–time" random set $\{(x,t) \in C \times [0,1] : x \in J_t\}$ rather than $\{t : J_t \neq \varnothing\}$.

PROOF OF THEOREM 1.4. The proof is similar to the proofs of the previous section, but with the spatial component taken into account. Let $I_t$, $J_t$ and $J_{n,t}$ be defined as in Theorem 1.2. As before, we have that $P(\exists t \in [0,1] : F_t \neq \varnothing) = 1$ if and only if $P(\exists t \in [0,1] : J_t \neq \varnothing) > 0$. Next, let $T_n := \{t : \bigcup_{k=1}^n I_{k,t} \neq C\}$ and note that $T_n$ is closed. Also, it is an elementary topology exercise (using the fact that the arcs are open and Brownian motion paths are continuous) left to the reader to check that if $t \in \bigcap_n T_n$, then $J_t \neq \varnothing$. Hence $P(A_n)$ bounded away from 0 implies that $P(\exists t \in [0,1] : J_t \neq \varnothing) > 0$.

Part (i). We will first show that if $\liminf_n n^3 u_n < \infty$, then

(16) $$\lim_{n \to \infty} P(\exists t \in [0,1] : J_{n,t} \neq \varnothing) = 0$$

[or equivalently $\lim_{n \to \infty} P(T_n \neq \varnothing) = 0$]. Fix $n$. For $i = 1, 2, \ldots, n^2$ and $j = 1, 2, \ldots, n$, put $A(i,j)$ for the event that for some $t \in [(i-1)/n^2, i/n^2]$ and some $x \in [(j-1)/n, j/n]$, $x \notin \bigcup_{k=1}^n I_{k,t}$. Then

$$P(\exists t \in [0,1] : J_{n,t} \neq \varnothing) \leq n^3 P(A(1,1)).$$

We have $A(1,1) \subset \bigcap_{k=1}^n B_k^c$ where $B_k$ is the event that $[0, 1/n] \subseteq I_{k,t}$ for every $t \in [0, 1/n^2]$. The event $B_k$ in turn contains the event that $I_{k,0} \supseteq [-M, M + 1/n]$, where, as in the above proof, $M = \max_{t \in [0,1/n^2]} |U_{k,t} - U_{k,0}|$. Letting $B_k'$ denote this last event, we have

$$P(B_k) \geq P(B_k') = \mathbb{E}[P(B_k'|M)] \geq \mathbb{E}\left[\ell_k - 2M - \frac{1}{n}\right] \geq \ell_k - \frac{5}{n}$$

and consequently

$$P(\exists t \in [0,1] : J_{n,t} \neq \varnothing) \leq n^3 \prod_{k=1}^n \left(1 - \ell_k + \frac{5}{n}\right)$$

$$\leq n^3 \prod_{k=1}^n (1 - \ell_k)\left(1 + \frac{10}{n}\right) \leq O(n^3 u_n).$$

Since the left-hand side is decreasing in $n$, (16) is established. For the case $\liminf_n n^3 u_n < \infty$, define $N_n$ as the number of $(i,j)$ for which $A(i,j)$ occurs. Then the above gives $\liminf_n \mathbb{E} N_n < \infty$. Letting

$$\tilde{T} := \{(x,t) \in C \times [0,1] : x \in J_t\},$$

we easily get $\tilde{T} \leq \liminf_n N_n$ and so by Fatou's lemma, we have that $\tilde{T}$ is a.s. finite. In particular, the set $\{t \in [0,1] : J_t \neq \varnothing\}$ is finite a.s. Again, Theorem 6.7 in [8] and (2.9) in [7] allow us to conclude that the latter set is empty a.s.

Part (ii). Let $\tilde{T}_n := \{(x,t) : x \in J_{n,t}\}$. Then the $\tilde{T}_n$'s are closed and $\bigcap_n \tilde{T}_n = \tilde{T}$; hence if all the $\tilde{T}_n$'s are nonempty, then so is $\tilde{T}$. Thus it suffices to show that $P(\tilde{T}_n \neq \varnothing)$ is bounded away from 0. Let $X_n$ be the two-dimensional Lebesgue measure of $\tilde{T}_n$. By Fubini's theorem, $\mathbb{E}X_n = u_n$ and so when $\mathbb{E}[X_n^2] \leq O(u_n^2)$ an application of the second moment method tells us that $\inf_n P(\tilde{T}_n \neq \varnothing) > 0$. We will now show that $\mathbb{E}[X_n^2] \leq O(u_n^2)$ if and only if $\sum_{n=1}^\infty e^{\ell_1 + \cdots + \ell_n}/n^4 < \infty$.

Fix $n$. Let $A_{t,x}$ be the event that $(t,x)$ is in $\tilde{T}_n$. Then, again by stationarity, $\mathbb{E}[X_n^2]$ is at most

$$2\int_C \int_0^1 P(A_{t,x} \cap A_{0,0})\,dt\,dx$$

and at least

$$\tfrac{1}{2}\int_C \int_0^{1/2} P(A_{t,x} \cap A_{0,0})\,dt\,dx.$$

Independence yields

$$P(A_{t,x} \cap A_{0,0}) = \prod_{k=1}^n P(A_{k,t,x} \cap A_{k,0,0}),$$

where $A_{k,t,x}$ is the event that $I_{k,t}$ does not contain $x$. Now $A_{k,t,x} \cap A_{k,0,0}$ is the event that $I_{k,0} \cap \{0, x - Z_t\} = \varnothing$, where $Z_t = U_{k,t} - U_{k,0}$ is the increment of $I_k$ in the time interval $[0,t]$, which is a normal random variable with mean 0 and variance $t$ projected onto $C$ as described in the Introduction. Hence

$$P(A_{k,t,x} \cap A_{k,0,0}) = \mathbb{E}[P(A_{k,t,x} \cap A_{k,0,0}|Z_t)] = 1 - 2\ell_k + \mathbb{E}[(\ell_k - |Z_t - x|)^+].$$

Inserting this into the product above gives, using the derivations in the proof of Theorem 1.2,

$$P(A_{t,x} \cap A_{0,0}) = \Theta(1)u_n^2 e^{\sum_{k=1}^n \mathbb{E}[(\ell_k - |Z_t - x|)^+]}.$$

Hence $\mathbb{E}[X_n^2] \leq O(u_n^2)$ if

$$\int_C \int_0^1 e^{\sum_{n=1}^\infty \mathbb{E}[(\ell_n - |Z_t - x|)^+]}\,dt\,dx < \infty \tag{17}$$

and only if

$$\int_C \int_0^{1/2} e^{\sum_{n=1}^\infty \mathbb{E}[(\ell_n - |Z_t - x|)^+]}\,dt\,dx < \infty. \tag{18}$$

The terms in the exponent are obviously bounded above by $\ell_n P(|Z_t - x| < \ell_n)$. For $t \geq 1/2$, this is clearly $O(1)\ell_n^2$ and so (17) holds if and only if (18) does.

Note next that for all positive $t$ and $x \in C$ we have, using Lemma 2.3 and basic facts of the normal distribution, that

$$\sum_{n:\ell_n^2 < t} \ell_n P(|Z_t - x| < \ell_n) = O(1) \sum_{n:\ell_n^2 < t} \frac{\ell_n^2}{\sqrt{t}} = O(1)$$

and

$$\sum_{n:\ell_n < 2|x|} \ell_n P(|Z_t - x| < \ell_n) = O(1) \sum_{n:\ell_n < 2|x|} \ell_n \frac{\ell_n}{|x|} = O(1).$$

To see that $P(|Z_t - x| < \ell_n) = O(1)\ell_n/|x|$, we use Lemma 2.3 together with the fact that the probability that a standard normal random variable takes values in $[y, y+d]$ with $y > 0$ is bounded above by $O(1)d/y$ and note that there is nothing to prove unless $\ell_n \leq |x|/2$.

It follows that (17) is equivalent to

$$\int_C \int_0^1 e^{\sum_{n:\ell_n \geq 2|x| \vee \sqrt{t}} \mathbb{E}[(\ell_n - |Z_t - x|)^+]} \, dt \, dx < \infty.$$

A lower bound for the exponent is given by

$$\sum_{n:\ell_n \geq 2|x| \vee \sqrt{t}} \mathbb{E}[(\ell_n - |Z_t - x|)^+]$$

$$= \sum_{n:\ell_n \geq 2|x| \vee \sqrt{t}} \mathbb{E}[(\ell_n - |Z_t - x|)^+ \mid |Z_t - x| < \ell_n] P(|Z_t - x| < \ell_n)$$

$$\geq \sum_{n:\ell_n \geq 2|x| \vee \sqrt{t}} (\ell_n - |x| - \sqrt{t}) P(|Z_t - x| < \ell_n)$$

$$= \sum_{n:\ell_n \geq 2|x| \vee \sqrt{t}} \ell_n P(|Z_t - x| < \ell_n) + O(1).$$

Thus our integral condition is equivalent to

$$\int_C \int_0^1 e^{\sum_{n:\ell_n \geq 2|x| \vee \sqrt{t}} \ell_n P(|Z_t - x| < \ell_n)} \, dt \, dx < \infty.$$

However, when $\ell_n \geq 2|x| \vee \sqrt{t}$,

$$P(|Z_t - x| < \ell_n) \geq P\left(|Z_t| \leq \frac{\ell_n}{2}\right) = 1 - O(1)\frac{\sqrt{t}}{\ell_n}.$$

Thus

$$\sum_{n:\ell_n \geq 2|x| \vee \sqrt{t}} \ell_n P(|Z_t - x| \geq \ell_n) = O(1) \sum_{n:\ell_n \geq \sqrt{t}} \sqrt{t} = O(1).$$

This shows that $\mathbb{E}[X_n^2] \leq O(u_n^2)$ if and only if

$$\int_C \int_0^1 e^{\sum_{n:\ell_n \geq 2|x| \vee \sqrt{t}} \ell_n} \, dt \, dx < \infty.$$

This clearly holds if and only if

$$\int_0^1 \int_0^1 e^{\sum_{n:\ell_n \geq x \vee \sqrt{t}} \ell_n} \, dt \, dx < \infty.$$

Now

$$\int_0^1 \int_0^1 e^{\sum_{n:\ell_n \geq x \vee \sqrt{t}} \ell_n} \, dt \, dx$$

$$= \int_0^1 \int_0^{x^2} e^{\sum_{n:\ell_n \geq x} \ell_n} \, dt \, dx + \int_0^1 \int_0^{\sqrt{t}} e^{\sum_{n:\ell_n \geq \sqrt{t}} \ell_n} \, dx \, dt$$

$$= \int_0^1 x^2 e^{\sum_{n:\ell_n \geq x} \ell_n} \, dx + \int_0^1 \sqrt{t} e^{\sum_{n:\ell_n \geq \sqrt{t}} \ell_n} \, dt = 3 \int_0^1 x^2 e^{\sum_{n:\ell_n \geq x} \ell_n} \, dx$$

where the last equality follows from the substitution $x = \sqrt{t}$. Via the substitution $u = x^3$ the last expression becomes

$$\int_0^1 e^{\sum_{n:\ell_n \geq u^{1/3}} \ell_n} \, du$$

and Lemma 2.1 now proves part (ii).

Part (iii)(a). We will more or less follow steps 1 and 2 in the previous HD arguments.

We first show that $\liminf_n n^\beta u_n < \infty$ implies

$$\mathrm{HD}(\tilde{T}) \leq \min\left\{3 - \beta, \frac{4-\beta}{2}\right\}.$$

Once this is done, the fact that $3 - x > 2 - \frac{x}{2}$ on $(0, 2)$ and the reverse holds on $(2, 3)$ and using Lemma 2.2, the upper bounds will be obtained. Consider now the union of the set of rectangles of the form $[(i-1)/n^2, i/n^2] \times [(j-1)/n, j/n]$ which contain a point $(t, x)$ with $x \notin \bigcup_{k=1}^n I_{k,t}$. This is a covering of $\tilde{T}$ with $N_n$ elements and from what we have seen in the proof of part (i) we can conclude that $\mathbb{E}[N_n] = O(1) n^3 u_n$. Since the elements of the covering have diameter of order $1/n$, we can conclude, as earlier, that $\mathrm{HD}(\tilde{T}) \leq 3 - \beta$. If we instead cover by $1/n^2 \times 1/n^2$ boxes, we get a covering $\tilde{T}$ with $N_n'$ elements of diameter of order $1/n^2$ with $\mathbb{E}[N_n'] = O(1) n^4 u_n$ and we can conclude that $\mathrm{HD}(\tilde{T}) \leq (4 - \beta)/2$ as well.

For the lower bound, assume first that $\beta_0 \in (0, 2)$. Then Theorem 1.2(iii) says that for each $x \in C$, $\mathrm{HD}(\{t \in [0, 1] : x \in F_t\}) = 1 - \beta_0/2$ a.s. By Fubini's theorem, we conclude that

$$\left\{x : \mathrm{HD}(\{t \in [0, 1] : x \in F_t\}) = 1 - \frac{\beta_0}{2}\right\}$$

has Lebesgue measure 1 a.s. Now, Theorem 7.7 in [14] (with $f$ there taken to be the projection onto $[0,1]$) allows us to conclude that $\mathrm{HD}(\{(t,x) : x \in F_t\}) \geq 2 - \beta_0/2$. (Theorem 7.7. says vaguely that any set in the square almost all of whose "slices" have HD $1 - \beta_0/2 > 0$ must have HD at least $2 - \beta_0/2$.)

For $\beta_0 \in [2,3)$, we argue differently. One follows the HD lower bound argument (suitably modified to a space–time situation) of Proposition A.16 of [12] mentioned earlier. In this case, one places a random measure on $\tilde{T}_n$ and obtains a uniform upper bound on the expected $3 - \beta_0 - \varepsilon$-energy. (Given a measure $m$ on $[0,1]^2$ and $\gamma > 0$, the $\gamma$-energy of $m$ is

$$\int_{[0,1]^2} \int_{[0,1]^2} |t-s|^{-\gamma} \, dm(t) \, dm(s).)$$

The random measure is of course Lebesgue measure restricted to $\tilde{T}_n$ and normalized by $u_n$. Using what was derived in part (ii), obtaining a uniform upper bound on the expected energy reduces to verifying the finiteness of

$$\int_0^1 \int_0^1 e^{\sum_{n : \ell_n \geq x \vee \sqrt{t}} \ell_n} \left(\frac{1}{x^2+t^2}\right)^{(3-\beta)/2} dt \, dx$$

under the assumption $\sum_{n=1}^\infty e^{\ell_1 + \ell_2 + \cdots + \ell_n}/n^{1+\beta} < \infty$. Breaking up the double integral as in the first part of the proof and checking that

$$\int_0^{x^2} \left(\frac{1}{x^2+t^2}\right)^{(3-\beta)/2} dt = \Theta(1) x^{\beta-1}$$

and

$$\int_0^{t^{1/2}} \left(\frac{1}{x^2+t^2}\right)^{(3-\beta)/2} dx = \Theta(1) t^{(\beta-2)/2},$$

it reduces to the finiteness of

$$\int_0^1 u^{\beta-1} e^{\sum_{n : \ell_n \geq u} \ell_n} \, du.$$

Another change of variables ($w = u^\beta$) together with Lemma 2.1 shows that finiteness of this integral is equivalent to the convergence of the given series.

Part (iii)(b): The $0 \leq \beta_0 < 2$ case follows from Theorem 1.2. The other cases follow from part (iii)(a) together with the fact that projections do not increase HD; see, for example, Theorem 7.5 in [14]. Alternatively, one can use a covering argument.

Part (iii)(c): The $0 \leq \beta_0 < 1$ case follows from Theorem 1.1. For the other cases, break the time interval into intervals of length $1/n^2$ and consider those intervals which contain a $t$ such that $\bigcup_{k=1}^n I_{k,t} \neq C$. If $N_n$ is the number of such intervals, we have from what was derived in part (i) that $E(N_n) = O(1) n^3 u_n$. This as before leads to the upper bound $(3-\beta_0)/2$ for the HD, as desired. $\square$

As for the type (I) case, the reason that we do not know if failure of the second moment method implies nonexistence of exceptional times is due to the fact that we cannot control the positions of the first $n$ arcs at the first time that the circle fails to be covered by them.

PROOF OF THEOREM 1.5. Part (i): Fix $n$ and partition $[0,1] \times C$ into boxes of size $1/n^\alpha \times 1/n$. For the given block $[0, 1/n^\alpha] \times [0, 1/n]$, we have, using arguments similar to those given earlier,

$$P\left(\exists (t,x) \in B : x \notin \bigcup_{k=1}^n I_{k,t}\right) \leq \prod_{k=1}^n P(\exists (t,x) \in B : x \notin I_{k,t})$$

$$\leq \prod_{k=1}^n \left(1 - \left(\ell_k - \frac{1}{n}\right) e^{-1/(n\ell_k)^\alpha}\right).$$

As before, one can show this is $O(1)u_n$. Since the number of blocks is of order $n^{1+\alpha}$, we get $P(\exists t \in [0,1] : F_t \neq \varnothing) = 0$ if $\liminf_n n^{1+\alpha} u_n = 0$ and proceed as earlier if $\liminf_n n^{1+\alpha} u_n \in (0, \infty)$.

Part (ii): We use the same notation as in Theorem 1.4. Using that argument (together with the analogous small modification given in Theorem 1.3 that dealt with the fact that certain time sets were not closed), proving the existence of exceptional times comes down to showing that $\mathbb{E}[X_n^2] \leq O(u_n^2)$. We now show that this holds if and only if the sum in the statement of the theorem is convergent.

By conditioning on whether arc $I_k$ has been updated by time $t$ or not we get

$$P(A_{k,t,x} \cap A_{k,0,0}) = (1 - e^{-t/\ell_k^\alpha})(1 - \ell_k)^2 + e^{-t/\ell_k^\alpha}(1 - 2\ell_k + (\ell_k - |x|)^+).$$

Hence

$$P(A_{t,x} \cap A_{0,0}) = u_n^2 \prod_{k=1}^n \left(1 - e^{-t/\ell_k^\alpha} + \frac{e^{-t/\ell_k^\alpha}(1 - 2\ell_k + (\ell_k - |x|)^+)}{(1-\ell_k)^2}\right)$$

$$= \Theta(1) u_n^2 \prod_{k=1}^n \left(1 + \frac{e^{-t/\ell_k^\alpha}(\ell_k - |x|)^+}{(1-\ell_k)^2}\right)$$

$$= \Theta(1) u_n^2 \prod_{k=1}^n (1 + e^{-t/\ell_k^\alpha}(\ell_k - |x|)^+)$$

$$= \Theta(1) u_n^2 e^{\sum_{k=1}^n e^{-t/\ell_k^\alpha}(\ell_k - |x|)^+}.$$

Thus $\mathbb{E}[X_n^2] \leq O(u_n^2)$ if and only if

$$\int_C \int_0^1 e^{\sum_{n=1}^\infty e^{-t/\ell_n^\alpha}(\ell_n - |x|)^+} \, dt \, dx < \infty$$

which clearly holds if and only if

$$\int_0^1 \int_0^1 e^{\sum_{n=1}^\infty e^{-t/\ell_n^\alpha}(\ell_n-x)^+}\, dt\, dx < \infty.$$

By a series of considerations analogous to what has been done in the earlier proofs we get

$$\sum_{n=1}^\infty e^{-t/\ell_n^\alpha}(\ell_n - x)^+ = \sum_{n:\,\ell_n^\alpha \geq t} e^{-t/\ell_n^\alpha}(\ell_n - x)^+ + O(1)$$
$$= \sum_{n:\,\ell_n^\alpha \geq t} (\ell_n - x)^+ + O(1) = \sum_{n:\,\ell_n^\alpha \geq x^\alpha \vee t} \ell_n + O(1).$$

Hence the given integral is finite if and only if

$$\int_0^1 \int_0^1 e^{\sum_{n:\,\ell_n^\alpha \geq x^\alpha \vee t} \ell_n}\, dt\, dx < \infty.$$

Copying the final parts of the proof of Theorem 1.4, we get

$$\int_0^1 \int_0^1 e^{\sum_{n:\,\ell_n^\alpha \geq x^\alpha \vee t} \ell_n}\, dt\, dx = \int_0^1 e^{\sum_{n:\,\ell_n^{\alpha+1} \geq u} \ell_n}\, du.$$

Now apply Lemma 2.1.

Part (iii): For the Hausdorff dimension upper bounds, we only sketch these. First, we have seen that $P(\tilde{T}_n \cap [0, 1/n^\alpha] \times [0, 1/n] \neq \varnothing) \leq O(u_n)$. Consider (a) and (a'). For $\alpha \geq 1$, we break up either into $(n^{-\alpha} \times n^{-1})$-boxes or $(n^{-\alpha} \times n^{-\alpha})$-boxes depending on whether $\beta_0$ is $\geq$ or $\leq \alpha$ and for $\alpha < 1$, we break up either into $(n^{-\alpha} \times n^{-1})$-boxes or $(n^{-1} \times n^{-1})$-boxes depending on whether $\beta_0$ is $\geq$ or $\leq 1$. This yields the upper bounds.

For (b), partition space into intervals of length $1/n$ and proceed in the same way. For (c), partition time into intervals of length $1/n^\alpha$ and proceed in the same way.

For the lower bounds, we follow the same arguments as in Theorem 1.4. First assume that $\alpha \geq 1$. If $\beta_0 \in (0, \alpha)$, we argue exactly as in the case $\beta_0 \in (0, 2)$ in Theorem 1.4 with 2 replaced by $\alpha$ throughout. If $\beta_0 \in [\alpha, \alpha+1)$, we argue exactly as in the case $\beta_0 \in [2, 3)$ in Theorem 1.4 where now things come down to verifying the finiteness of

$$\int_0^1 \int_0^1 e^{\sum_{n:\,\ell_n^\alpha \geq x^\alpha \vee t} \ell_n} \left(\frac{1}{x^2 + t^2}\right)^{(1+\alpha-\beta)/2} dt\, dx$$

under the assumption $\sum_{n=1}^\infty e^{\ell_1 + \ell_2 + \cdots + \ell_n}/n^{1+\beta} < \infty$. Now assume that $\alpha < 1$. If $\beta_0 \in (0, 1)$, we argue analogously but a theorem of Kahane replaces our use of Theorem 1.2(iii). Theorem 4 in Section 11.8 of [11] together with Lemma 2.2 tells us that for each $t \in C$, $\text{HD}(\{x \in C : x \in F_t\}) = 1 - \beta_0$ a.s. An application of Fubini's theorem and Theorem 7.7 in [14] exactly as in the case $\beta_0 \in (0,2)$ in Theorem 1.4 allows us to conclude $\mathrm{HD}(\{(t,x): x \in F_t\}) \geq 2 - \beta_0$, as desired. For $\beta_0 \in [1, 1+\alpha)$, we argue as in the case $\beta_0 \in [2,3)$ in Theorem 1.4 where now things come down to verifying the finiteness of

$$\int_0^1 \int_0^1 e^{\sum_{n:\ell_n^\alpha \geq x^\alpha \vee t} \ell_n} \left(\frac{1}{x^2 + t^2}\right)^{(1+\alpha-\beta)/(2\alpha)} dt\, dx$$

under the assumption $\sum_{n=1}^\infty e^{\ell_1 + \ell_2 + \cdots + \ell_n}/n^{1+\beta} < \infty$. This is done in more or less the same way with a few easy needed modifications. $\square$

PROOF OF THEOREM 1.6. Part (i): Fix $n$ and partition $[0,1] \times C$ into boxes of size $1/\log n \times 1/n$. For the given block $B$, we have, by using the same arguments as in Theorem 1.5,

$$P\left(\exists (t,x) \in B : x \notin \bigcup_{k=1}^n I_{k,t}\right) \leq \prod_{k=1}^n P(\exists (t,x) \in B : x \notin I_{k,t})$$
$$\leq \prod_{k=1}^n \left(1 - \left(\ell_k - \frac{1}{n}\right) e^{-1/\log n}\right).$$

One can again show this is $O(1) u_n$. Since the number of blocks is of order $n \log n$, we get $P(\exists t \in [0,1]: F_t \neq \varnothing) = 0$ if $\liminf_n n(\log n) u_n = 0$ and proceed as earlier if $\liminf_n n(\log n) u_n \in (0, \infty)$.

Part (ii): The difference between the situation here and that of the previous proof is that Lemma 2.1 does not work for $\beta = 0$. Therefore the analysis will be slightly different even though the ideas are the same. By repeating the beginning of the previous proof, it follows that existence of exceptional times is implied by

$$\int_0^1 \int_0^1 e^{e^{-t} \sum_{n=1}^\infty (\ell_n - x)^+} dt\, dx < \infty$$

which in turn is, of course, equivalent to

$$\int_0^b \int_0^b e^{e^{-t} \sum_{n=1}^\infty (\ell_n - x)^+} dt\, dx < \infty$$

for any $b > 0$. The rest of the proof will consist of showing analytically that this is equivalent to the convergence of the series in the statement of the theorem. Since $1 - t < e^{-t} < 1 - t/2$ [for $t \in (0,1)$], it suffices to show that for any $\kappa \in [1/2, 1]$,

$$\int_0^b \int_0^b e^{(1-\kappa t) \sum_{n=1}^\infty (\ell_n - x)^+} dt\, dx < \infty$$

is equivalent to convergence of the series. Let $g(x) := \sum_{n=1}^\infty (\ell_n - x)^+$ and choose $b$ so that $g(x) \geq 2$ on $[0, b]$. Now

$$\int_0^b e^{(1-\kappa t) \sum_{n=1}^\infty (\ell_n - x)^+} dt = e^{g(x)} \left[-\frac{e^{-\kappa t g(x)}}{\kappa g(x)}\right]_0^b = e^{g(x)} \frac{1 - e^{-b\kappa g(x)}}{\kappa g(x)}.$$

The last expression is $\Theta(1)e^{g(x)}/(g(x))$. Thus

$$\int_0^b \int_0^b e^{(1-\kappa t)\sum_{n:\ell_n \geq x} \ell_n} \, dt \, dx = \Theta(1) \int_0^b \frac{e^{g(x)}}{g(x)} \, dx = \Theta(1) \int_0^b e^{f(x)} \, dx,$$

where $f(x) := g(x) - \log g(x)$. We again use Lemma 11.4.1 of [11] which was stated in the proof of Lemma 2.1. We have that

$$f'(x) = g'(x)\left(1 - \frac{1}{g(x)}\right) = -n\left(1 - \frac{1}{\ell_1 + \ell_2 + \cdots + \ell_n}\right), \qquad \ell_{n+1} < x < \ell_n.$$

By choice of $b$, $f'$ is negative and increasing. Hence $f$ is convex and decreasing and we may apply the lemma. Since $f'(x) = \Theta(n)$, $\ell_{n+1} < x < \ell_n$, $f'$ makes jumps of size $1 + o(1)$ when $x = \ell_n$ and

$$f(\ell_n) = \ell_1 + \ell_2 + \cdots + \ell_n - n\ell_n - \log(\ell_1 + \ell_2 + \cdots + \ell_n - n\ell_n)$$
$$= \ell_1 + \ell_2 + \cdots + \ell_n - \log \log n + O(1),$$

the lemma gives us that

$$\int_0^b e^{f(x)} \, dx = \sum_{n=1}^\infty \frac{e^{\ell_1 + \ell_2 + \cdots + \ell_n}}{n^2 \log n},$$

as desired. $\square$

**4. Further questions.** In this section, we list a number of questions and problems that remain:

1. When $\ell_n = 1/n$ and $\alpha = 0$, are there exceptional times in the Poisson model?

2. Show that the inequalities in Theorem 1.4(iii) are equalities.

3. If $c < 1$ and $\ell_n = c/n$, then we know that $P(F = \varnothing) = 0$. Is it also the case that

$$P(\exists t \in [0,1] : F_t = \varnothing) = 0?$$

Does this depend on the value of $c$?

This is analogous to the dynamical percolation question of whether, when we do percolate for ordinary percolation, there are exceptional times at which percolation does *not* occur. For dynamical percolation, this question is much less understood than the reverse question where one does not percolate for ordinary percolation but asks if there are exceptional times at which percolation *does* occur.

4. Given subsets of the time interval, determine when they contain exceptional times of various types.

MATHEMATICAL SCIENCES
CHALMERS UNIVERSITY OF TECHNOLOGY
AND
MATHEMATICAL SCIENCES
GÖTEBORG UNIVERSITY
SE-41296 GOTHENBURG
SWEDEN
E-MAIL: jonasson@math.chalmers.se
steif@math.chalmers.se